# ВИКОРИСТАННЯ ДИФЕРЕНЦІАЛЬНИХ РІВНЯНЬ ЗІ ЗМІННИМИ КОЕФІЦІЄНТАМИ ДЛЯ ОПИСУ РУХІВ НЕЛІНІЙНИХ ЕЛЕКТРОМЕХАНІЧНИХ СИСТЕМ

**Волянський Р.С.**, к.т.н., доц.
*КПІ ім. Ігоря Сікорського, кафедра автоматизації електромеханічних систем та електроприводу*

**Вступ.** Внаслідок процесів, які відбуваються при функціонуванні сучасних електромеханічних систем (ЕМС), ці системи з точки зору теорії динамічних систем можна розглядати як складні нелінійні динамічні системи [1-3]. Рух таких систем повністю визначається зовнішніми впливами, які діють на ЕМС, їх параметрами та початковими умовами роботи [4-6]. Причому в деяких випадках, які обумовлені певними значеннями параметрів та початкових умов роботи ЕМС, їх динаміка перестає бути регулярною [7-9] та набуває хаотичного характеру [10-13].

Вищезазначені фактори ускладнюють дослідження електромеханічних систем та, в загальному випадку, унеможливлюють використання класичних методів аналізу динаміки ЕМС [14-15], оскільки останні нехтують особливостями нелінійних систем [16-18] та описують їх динаміку за допомогою звичайних лінійних диференційних рівнянь зі сталими коефіцієнтами. В той же час в теорії керування напрацьовано велику кількість різноманітних методів і підходів для аналізу стійкості [19-20] та синтезу траєкторій руху саме на основі лінійних диференційних рівнянь [21-22].

Тому виникає важлива задача створення математичних моделей нелінійних ЕМС, які врахують особливості їх руху, але мають вигляд аналогічний до лінійних моделей.

**Метою роботи** є розробка методики опису руху нелінійних електромеханічних систем за допомогою лінійних диференційних рівнянь зі змінними коефіцієнтами, які залежать від координат досліджуваної ЕМС.

**Матеріали та результати досліджень.** В нашій роботі ми будемо розглядати рух узагальненої електромеханічної системи, яка описується диференційними рівняннями виду:

$$\frac{d}{dt}\mathbf{y} = \mathbf{f}(\mathbf{y}, \mathbf{u}, t), \qquad (1)$$

тут: **f** – матриця-функція, **y**, **u** – вектори змінних стану та зовнішні впливи, t – час

Наявність матриці-функції f() робить цю систему нелінійною у випадках, коли кожен компонент $f_i()$ матриці-функції **f**() має вигляд складніший за зважену за деякими коефіцієнтами лінійну комбінацію змінних стану розглядаємої системи.

Як вже зазначалося, класичний підхід до аналізу нелінійних систем передбачає здійснювати розкладання кожного компоненти функції **f** в ряд Тейлора в околиці робочої точки [21]. При цьому система (1) перетворюється до вигляду:

$$\frac{d}{dt}\mathbf{y} = \mathbf{By} + \mathbf{mU} + \mathbf{g}(\mathbf{y},\mathbf{U},t), \quad (2)$$

де матриці коефіцієнтів **B** и **m** визначаються наступними залежностями:

$$\mathbf{B} = \left.\frac{\partial \mathbf{f}}{\partial \mathbf{y}}\right|_{\mathbf{y}=\mathbf{y}^*, \mathbf{u}=\mathbf{u}^*}; \quad \mathbf{m} = \left.\frac{\partial \mathbf{f}}{\partial \mathbf{u}}\right|_{\mathbf{y}=\mathbf{y}^*, \mathbf{u}=\mathbf{u}^*}, \quad (3)$$

а $\mathbf{g}(\mathbf{y},\mathbf{U},t)$ – є матрицею-функцією, компоненти якої є членам ряду Тейлора другого і вищих порядків, $\mathbf{y}^*$ та $\mathbf{u}^*$ є векторами, які визначають положення робочої точки розглядаємої системи у фазовому просторі та керуючі сигнали, які подаються на систему в цій точці.

Оскільки останній доданок у рівняннях (2) визначається поліноміальною залежністю, то при задачах аналізу та синтезу руху цим доданком нехтують

$$\frac{d}{dt}\mathbf{y} = \mathbf{By} + \mathbf{mU}. \quad (4)$$

Очевидно, що якщо об'єкт (1) описується істотно нелінійними функціями $f_i()$, то при виконанні подібної процедури відбувається значна втрата інформації про розглядаєму систему. Таким чином, (4) є значним спрощенням (1) і отримані шляхом використання рівняння (4) результати взагалі-то є справедливими лише для робочої точки системи (1). При цьому не гарантується досягнення цієї точки на траєкторіях руху (4) та збереження стійкості руху в цілому.

Представлення нелінійної функції у виді перемінних коефіцієнтів, що визначаються на основі методу гармонійної лінеаризації нелінійної залежності

$$q(A,\omega_c) = \int_0^{2\pi}\left[f_i\big(A\sin(\omega_c t),\omega_c t\big)\sin(\omega_c t)\right]d\omega_c t;$$
$$q'(A,\omega_c) = \int_0^{2\pi}\left[f_i\big(A\sin(\omega_c t),\omega_c t\big)\cos(\omega_c t)\right]d\omega_c t, \quad (5)$$

де $A$ и $\omega_c$ – амплітуда і частота еквівалентного лінеаризуючого впливу відповідно, приводить до перетворення системи рівнянь (1) до наступного вигляду:

$$\frac{d}{dt}\mathbf{y} = \frac{\mathbf{q}(\mathbf{y},\omega_c)\cdot\boldsymbol{\eta} + \mathbf{q}(U,\omega_c)\cdot\mathbf{U}}{\mathbf{q}'(\mathbf{y},\omega_c)+1}. \qquad (6)$$

Такий підхід дає досить точні результати для простих нелінійних функцій, що не містять петльових залежностей. Однак представлення узагальненої функції $\mathbf{f}(\mathbf{y},\mathbf{U},t)$ у виді $\mathbf{f}(A\sin(\omega t),\omega t)$ вимагає виконання ряду додаткових обчислень, які складно формалізуються і являють собою самостійну досить складну задачу. До того ж, метод гармонійної лінеаризації є наближеним через припущення про синусоїдальну зміну координат розглядаємої системи під дією зовнішнього синусоїдального впливу. Для багатьох нелінійних систем таке наближення є досить грубим оскільки формуємі траєкторії руху в наслідок наявності нелінійних залежностей можуть суттєво відрізнятися від синусоїдальних. Через цей факт рівняння (6) є черговим, більш точним у порівнянні з використанням системи (4), наближенням до рівнянь (1). Проте точної відповідності цих рівнянь немає.

Встановлювати точну відповідність між вихідними та перетвореними рівняннями ми пропонуємо використовуючи поняття псевдо-аффінних динамічних систем [21]. З використанням цього підходу матриці коефіцієнтів **B** и **m** для динамічної системи (1) будемо шукати в такий спосіб.

Помножимо і розділимо праву частину системи (1) на суму лінійних комбінацій керуючих впливів і координат об'єкта:

$$\frac{d}{dt}y_i = f_i(\mathbf{y},\mathbf{U},t)\frac{\sum_{i=1}^{n}a_{ik}y_i + \sum_{i=1}^{l}c_{ik}U_i}{\sum_{i=1}^{n}a_{ik}y_i + \sum_{i=1}^{l}c_{ik}U_i} =$$

$$= \left(\frac{f_i(\mathbf{y},\mathbf{U},t)}{\sum_{i=1}^{n}a_{ik}y_i + \sum_{i=1}^{l}c_{ik}U_i}\right)\sum_{i=1}^{n}a_{ik}y_i + \left(\frac{f_i(\mathbf{y},\mathbf{U},t)}{\sum_{i=1}^{n}a_{ik}y_i + \sum_{i=1}^{l}c_{ik}U_i}\right)\sum_{i=1}^{l}c_{ik}U_i. \qquad (7)$$

Коефіцієнти $a_{ik}, c_{ik}$ є довільними і повинні вибиратися в такий спосіб що б знаменник (7) не дорівнював нулю в усьому діапазоні можливих значень координат та керуючих впливів об'єкта (1).

Якщо функція $f_i(\mathbf{y},\mathbf{U},t)$ має строгий математичний опис, то подальший порядок перетворень буде наступним:

Введемо наступне позначення:

$$K_i(\mathbf{y},\mathbf{U},t) = \frac{f_i(\mathbf{y},\mathbf{U},t)}{\sum_{i=1}^{n} a_{ik} y_i + \sum_{i=1}^{l} c_{ik} U_i} \qquad (8)$$

та, підставивши (8) у (7), одержимо

$$\frac{d}{dt} y_i = K_i(\mathbf{y},\mathbf{U},t) \sum_{i=1}^{n} a_{ik} y_i + K_i(\mathbf{y},\mathbf{U},t) \sum_{i=1}^{l} c_{ik} U_i, \qquad (9)$$

або

$$\frac{d}{dt} y_i = \sum_{i=1}^{n} b_{ik} y_i + \sum_{i=1}^{l} m_{ik} U_i, \qquad (10)$$

де

$$b_{ik} = a_{ik} \cdot K_i(\mathbf{y},\mathbf{U},t), \quad m_{ik} = c_{ik} \cdot K_i(\mathbf{y},\mathbf{U},t). \qquad (11)$$

У випадку, якщо функція $f_i(\mathbf{y},\mathbf{U},t)$ не має суворого математичного опису, такий опис є занадто складний, її можна замінити, відповідно до системи (1) наступним чином

$$\frac{d}{dt} y_i = f_i(\mathbf{y},\mathbf{U},t) \frac{\sum_{i=1}^{n} a_{ik} y_i + \sum_{i=1}^{l} c_{ik} U_i}{\sum_{i=1}^{n} a_{ik} y_i + \sum_{i=1}^{l} c_{ik} U_i} =$$
$$= \left( \frac{\frac{d}{dt} y_i}{\sum_{i=1}^{n} a_{ik} y_i + \sum_{i=1}^{l} c_{ik} U_i} \right) \sum_{i=1}^{n} a_{ik} y_i + \left( \frac{\frac{d}{dt} y_i}{\sum_{i=1}^{n} a_{ik} y_i + \sum_{i=1}^{l} c_{ik} U_i} \right) \sum_{i=1}^{l} c_{ik} U_i. \qquad (12)$$

Позначивши вираз в дужках наступним чином

$$K_i(\mathbf{y},\mathbf{U},t) = \frac{\frac{d}{dt} y_i}{\sum_{i=1}^{n} a_{ik} y_i + \sum_{i=1}^{l} c_{ik} U_i}, \qquad (13)$$

та визначивши коефіцієнти (11) одержимо систему рівнянь, яка є аналогічною до системи (10).

**Висновки.** Таким чином, у результаті виконання перетворень (7) – (9), (11)-(13) система нелінійних диференціальних рівнянь (1) приведена до тотожної системи (10) зі змінними коефіцієнтами (1). Причому ці перетворення легко формалізуються і можуть бути виконані за допомогою сучасних програмних засобів. Завдяки тому, що система (10) має вид аналогічний (4), аналіз та синтез траєкторій руху розглядаємої системи може здійснюватися будь-яким відомим методом.

### Перелік посилань

1. Садовой О.В. Синтез оптимальної системи керування з іраціональною активаційною функцією // Вісник НТУ «ХПІ» «Проблеми автоматизованого електроприводу Теорія та практика)». – Харків: НТУ «ХПИ». – 2010.– Вип. 28. – С.49 – 51 (на рос.)
2. Садовой О.В. Синтез оптимальної системи керування з нелінійною активаційною функцією. Електротехнічні і комп'ютерні системи. Київ: Техніка. 2014, № 15 (91), С. 69 – 71 (на рос.)
3. Цабенко М.В. Математичні моделі усунення помпажа в відцентровому компресорі // Вісник КДУ ім. М. Остоградського. - 2010. - Випуск 4(63). - С. 167 - 169. (на рос.)
4. Назарова Е.С. Система оптимального керування натягом прокатуваної полоси стану холодної прокатки // Збірка наукових праць Донбаського державного технічного університету. –2011. – Вип. 34. – С.122-130 (на рос.)
5. Волянський Р.С. Система керування слідкуючим електроприводом із ковзним режимом 2-го порядка, Вісник Кременчуцького державного університету імені Михайла Остроградського. – Кременчук: КДУ, 2010. – Вип. 4/2010 (63) частина 3. – с.11-14 (на рос.)
6. Sokhina, Y Active Suspension Control System. In: IEEE International Conference on Modern Electrical and Energy Systems (MEES), Kremenchuk, Ukraine, pp. 10-13 (2019)
7. Волянський Р. С. Синтез керуючого впливу для електромеханічних об'єктів з невідомими параметрами // Вісник НТУ «ХПІ». — Харьков: ХПИ. — 2015. — № 12 (1121). — С. 60—63 (на рос.)
8. Pranolo A, Parallel Mathematical Models of Dynamic Objects, Int. J. Adv. Intell. Informs. 4(2) (2018) 120–131.
9. I. Shramko, "Anti-swing Control System for the One Class of Underactuated Dynamic Objects," 2020 IEEE Problems of Automated Electrodrive. Theory and Practice (PAEP), Kremenchuk, Ukraine, 2020, pp. 1-4.
10. Kluev, O., "Chaotic Time-variant Dynamical System," 2020 IEEE 15th International Conference on Advanced Trends in Radioelectronics, Telecommunications and Computer Engineering (TCSET), 2020, pp. 606-609.
11. R. Voliansky, "Chua's circuits interval synchronization," 2017 4th International Scientific-Practical Conference Problems of Infocommunications. Science and Technology (PIC S&T), Kharkov, Ukraine, 2017, pp. 439-443.
12. O. Sinkevych, "Sliding Mode Control for DC Generator with Uncertain Load," 2020 IEEE 15th International Conference on Advanced Trends in Radioelectronics, Telecommunications and Computer Engineering (TCSET), Lviv-Slavske, Ukraine, 2020, pp. 313-316.

13. Voliansky, R. (2017). Transformation of the generalized chaotic system into canonical form. International Journal of Advances in Intelligent Informatics, 3(3), 117–124.
14. Sadovoi A. V., "The transformation of linear dynamical object's equation to Brunovsky canonical form," 2017 IEEE 4th International Conference Actual Problems of Unmanned Aerial Vehicles Developments (APUAVD), Kiev, Ukraine, 2017, pp. 196-199.
15. Садовой А. В. Дослідження частотних характеристик динамічних ланок із похідним дробових порядків / Р. С. Волянский, // Збірка наукових праць ДДТУ (технічні науки). Тематичний випуск „Проблеми автоматизованого електроприводу. Теорія і практика". – 2007. – С. 82–85. (на рос.)
16. O. Sinkevych, "Construction of Parallel Piecewise-Linear Interval Models for Nonlinear Dynamical Objects," 2019 9th International Conference on Advanced Computer Information Technologies (ACIT), 2019, pp. 97-100.
17. Y. Sokhina, "Sliding Mode Interval Controller for the Mobile Robot," 2019 XIth International Scientific and Practical Conference on Electronics and Information Technologies (ELIT), Lviv, Ukraine, 2019, pp. 76-81.
18. N. Volianska, "Interval model of the piezoelectric drive," 2018 14th International Conference on Advanced Trends in Radioelecrtronics, Telecommunications and Computer Engineering (TCSET), Lviv-Slavske, Ukraine, 2018, pp. 1-6.
19. O. Sinkevych, "Root Methods for Dynamic Analysis of the One Class Chaotic Systems," 2019 IEEE 14th International Conference on Computer Sciences and Information Technologies (CSIT), 2019, pp. 117-121.
20. O. Sadovoi, "Defining of Lyapunov Functions for the Generalized Nonlinear Object," 2018 IEEE 5th International Conference on Methods and Systems of Navigation and Motion Control (MSNMC), Kiev, Ukraine, 2018, pp. 222-228.
21. Садовой О. В. Лінеаризація зворотними зв'язками рівнянь динаміки узагальненого електромеханічного об'єкта з диференціальним рівнянням спостережуваності // Вісник Нац. техн. ун-ту "ХПІ": зб. науч. пр. Темат. вип. : Механіко-технологічні системи та комплекси. – Харків: НТУ "ХПІ". – 2014. – № 60 (1102). – С. 52-57. (на рос.)
22. N. Krasnoshapka, O. Statsenko, "Electromechanical System Motion Control in Direct and Inverse Time," 2022 IEEE 3rd KhPI Week on Advanced Technology (KhPIWeek), Kharkiv, Ukraine, 2022, pp. 1-6.